\documentclass[12pt]{article}

\usepackage{times, amsfonts, amstext, amsmath, amssymb, amsthm, titlesec}

\usepackage[authoryear]{natbib}

\usepackage{setspace}
\usepackage[breaklinks, colorlinks=true,allcolors=blue]{hyperref}
\usepackage{etoolbox}

\makeatletter

\patchcmd{\NAT@citex}
{\@citea\NAT@hyper@{%
		\NAT@nmfmt{\NAT@nm}%
		\hyper@natlinkbreak{\NAT@aysep\NAT@spacechar}{\@citeb\@extra@b@citeb}%
		\NAT@date}}
{\@citea\NAT@nmfmt{\NAT@nm}%
	\NAT@aysep\NAT@spacechar\NAT@hyper@{\NAT@date}}{}{}

\patchcmd{\NAT@citex}
{\@citea\NAT@hyper@{%
		\NAT@nmfmt{\NAT@nm}%
		\hyper@natlinkbreak{\NAT@spacechar\NAT@@open\if*#1*\else#1\NAT@spacechar\fi}%
		{\@citeb\@extra@b@citeb}%
		\NAT@date}}
{\@citea\NAT@nmfmt{\NAT@nm}%
	\NAT@spacechar\NAT@@open\if*#1*\else#1\NAT@spacechar\fi\NAT@hyper@{\NAT@date}}
{}{}
\makeatother

\usepackage{tabularx, siunitx, booktabs}
\usepackage[labelfont=bf, labelsep=period, singlelinecheck=false]{caption}

\long\def\symbolfootnote[#1]#2{\begingroup%
\def\thefootnote{\fnsymbol{footnote}}\footnote[#1]{#2}\endgroup}

\titleformat{\section}{\large\bfseries}{\thesection.}{.5em}{}
\titlespacing*{\section}{0pt}{*3}{*2}
\titleformat{\subsection}{\normalfont\bfseries}{\thesubsection.}{.5em}{}
\titlespacing*{\subsection} {0pt}{*3}{*2}
\titleformat{\subsubsection}{\normalfont\bfseries}{\thesubsubsection.}{.5em}{}
\titlespacing*{\subsubsection} {0pt}{*3}{*2}

\theoremstyle{plain}
\newtheorem{theorem}{Theorem}[section]
\newtheorem{lemma}{Lemma}[section]

\numberwithin{equation}{section}

\begin{document}

\title{\textbf{\LARGE Sequential estimation in the group testing problem}}

\date{}

\begingroup
\let\center\flushleft
\let\endcenter\endflushleft
\maketitle
\endgroup

\author{
\vskip -1cm
\noindent
{\large Gregory Haber\textsuperscript{a}, Yaakov Malinovsky\textsuperscript{a}, and Paul S. Albert\textsuperscript{b}}

\noindent \textsuperscript{a}Department of Mathematics and Statistics, University of Maryland, Baltimore County, Baltimore, Maryland, USA; \textsuperscript{b}Biostatistics Branch, Division of Cancer Epidemiology and Genetics, National Cancer Institute, Rockville, Maryland, USA}

\symbolfootnote[0]{\normalsize \textbf{CONTACT} Gregory Haber \href{mailto:ghaber1@umbc.edu}{ghaber1@umbc.edu}
Department of Mathematics and Statistics, University of Maryland, Baltimore County, 1000 Hilltop Circle, Baltimore,
MD 21250, USA}

{\small \noindent\textbf{Abstract:} Estimation using pooled sampling has long been an area of interest in the group testing literature. Such research has focused primarily on the assumed use of fixed sampling plans (i), although some recent papers have suggested alternative sequential designs that sample until a predetermined number of positive tests (ii). One major consideration, including in the new work on sequential plans, is the construction of debiased estimators which either reduce or keep the mean square error from inflating. Whether, however, under the above or other sampling designs unbiased estimation is in fact possible has yet to be established in the literature. In this paper, we introduce a design which samples until a fixed number of negatives (iii), and show that an unbiased estimator exists under this model, while unbiased estimation is not possible for either of the preceding designs (i) and (ii). We present new estimators under the different sampling plans that are either unbiased or that have reduced bias relative to those already in use as well as generally improve on the mean square error. Numerical studies are done in order to compare designs in terms of bias and mean square error under practical situations with small and medium sample sizes.}
\\ \\
{\small \noindent\textbf{Keywords:} Binomial sampling plans; group testing; sequential estimation; unbiased estimation.}
\\ \\
{\small \noindent\textbf{Subject Classifications:} 62F10; 62K05; 62L12.}

\section{Introduction}
\label{s:intro}
Group testing was first introduced in the statistical literature by \citet{dorfman1943} as an efficient blood screening technique, and has since developed into two nearly distinct research areas: classification and estimation. The main feature in both fields is that, for suitable ranges of the underlying population parameter, group testing offers the potential for significant savings in efficiency (in terms of reduced mean squared error (MSE)), the necessary number of trials, or both.

The basic problem addressed by Dorfman was as follows, suppose there is a large population, with each individual possessing some trait with probability $p$ and it is necessary to determine each person carrying said trait. The standard method of testing each individual would be prohibitive in terms of both time and cost, so Dorfman proposed to first group the units being tested into pools of size $k$. If this pool is tested and found negative, we can assume all individuals within it to be free of the trait. If it is positive, however, at least one carries the trait and each unit from that group would then be tested without pooling. In total, this means that for any given pool either $1$ or $k + 1$ tests are required as opposed to the $k$ needed under the standard method. Intuitively, if $p$ is small enough so that sufficiently many of the groups are negative, this can lead to a real savings in the total number of tests.

The primary question when using this method, known as Dorfman's two stage procedure, is to find a group size $k$ which minimizes the number of tests required for a given (assumed known) value of $p$. While \citet{dorfman1943} gave only numerical results, this problem and method became the framework for what is known as the classification, or identification, problem in group testing. Subsequent research has generalized to a wide array of designs and applications \citep[see, for example,][]{sobel1959, sobel1966, graff1972, hwang1976, johnson1990, barlev1995, barlev2005, ahl2012, mcm2012, tatsuoka, malinovsky2015, malinovsky2016}.

A second area in which group testing methods have developed is the estimation of a Bernoulli parameter $p$ using pooled samples. Many of the early applications of this idea took place prior to the concept being introduced in the statistical literature, with a classic example being the rate of disease transmission from an insect to a given plant  \citep[for specific cases, see][]{gibbs1960, thompson1962}. Such examples were motivated by a scarcity of available specimens for testing with the main idea being that if only a few plants were able to be tested, but many insects were available, it would be possible to capture the information of disease transmission from a larger number of insects if a pool of size $k$, as opposed to a single unit,  were exposed to the plant of interest.

Early studies of the statistical properties of this procedure were carried out by \citet{gibbs1960}, who addressed the bias of the standard MLE estimator, as well as \citet{thompson1962} and \citet{chiang1962}. Much subsequent research related to this procedure dealt with design issues, most importantly the selection of the group size $k$, based on two primary aims: reducing the probability of achieving all positive or negative groups \citep{kerr1971, griffiths1972} and the use of prior information to choose $k$ minimizing the MSE \citep{swallow1985, ho1994}.

Other work related to estimation has considered optimal resampling techniques \citep{sobel1975, chen1990}, finite population cases \citep{bhat1979, theobold2014}, testing with misclassification \citep{tu1995, liu2012}, and how to address multiple group sizes across samples \citep{walter1980, le1981, hepworth1996, dres2015, santos2016}. More recent extensions have included multinomial populations \citep{ho2000, tebbs2013, warasi2016} and regression techniques for non-homogeneous binary populations \citep{xie2001, bilder2009, chen2009}.

The use of sequential designs for group testing estimation was first made, to our knowledge, in \citet{kerr1971} who suggested that if a fixed size sample is taken and all groups are found to be positive or negative, new samples of predetermined size be collected iteratively until a desirable result is obtained. Later, \citet{ku2006} suggested the use of inverse binomial sampling (drawing until a certain number of positives) as a means of efficiently estimating $p$ in certain cases, particularly when the prevalence is small and a fixed design would require an unreasonably large sample size. \citet{pritchard2011a} generalized the same model to the case of unequal group sizes and provided debiased estimators for the uniform group-size model. Similarly, \citet{hepworth2013} generalized his own debiased estimator from the fixed design  \citep[see][]{hepworth2009} to the inverse binomial model. In addition, \citet{hepworth2013} proposed confidence interval estimators for the negative binomial model. Other work, all assuming the inverse binomial model, includes Bayesian estimation \citep{pritchard2011b} and regression methods \citep{mont2015}.

To carry out estimation, the maximum likelihood estimator (MLE) has been heavily used, justified by its simplicity and good large sample properties such as consistency and asymptotic normality. However, as mentioned above, the MLE is a biased estimator, quite heavily for some values of $p$ \citep{gibbs1960, thompson1962, swallow1985}. This issue arises frequently in applications, where the number of pools sampled is often too small to rely on the above properties. Fairly recent examples of this can be found in  \citet{williams2010} who analyzed two data sets related to disease prevalence in wild fish populations with 12 and 27 pools, and \citet{galetto2014} who examined the acquisition of Flavescence dor{\'e}e for different grapvine species from the leafhopper {\it Scaphoideus titanus} using samples with 21 and 24 pools.

To address this, alternative estimators have been proposed in attempts to reduce this bias such as in \citet{burrows1987} and \citet{hepworth2009}, as well as the above mentioned examples in the sequential case. In particular, Burrows was able to show empirically that his estimator not only improves on the bias, but, perhaps more importantly, yields a smaller MSE than the MLE for all values of $p$ considered ($p \leq 0.5)$.

Because so much emphasis is placed on the bias, it is important to note that, when standard binomial sampling is used, no unbiased estimator of $p$ exists. Although not mentioned often, this fact was stated by \citet{hall1963} in the context of group testing and proved for the general binomial case in \citet[p. 100]{lehmann1998}.
The question remains, however, as to whether an unbiased estimator exists when the broader class of binomial sampling plans (defined below in Section \ref{s:bsp}) are considered.

In this paper, we introduce a design which samples until a given number of negatives and show, based on results in \citet{degroot1959}, that unbiased estimation is possible for such a plan. It is also shown that unbiased estimation is, in fact, impossible under the previously mentioned designs found in the literature. Additionally, we extend the idea of Burrows for the fixed sample design to the case of inverse binomial sampling to provide reasonable debiased estimators for the sequential case, particularly for the model under which unbiased estimation is not feasible. Numerical studies are then done in order to compare designs/estimators in terms of bias and MSE under practical situations.

\section{Binomial sampling plans}
\label{s:bsp}
Before describing the models used here, we first give a definition of binomial sampling plans. A more comprehensive treatment, as well as several results on estimation under such models, can be found in \citet{girshick1946}.

For our purposes a general binomial sampling plan $S$ can be defined as a subset of the non-negative integer valued coordinates in the $xy-$plane determined by a set of boundary points, $B_S$, at which the sampling terminates. All plans begin at the origin and, until a point $\gamma = (X(\gamma), Y(\gamma)) \in B_S$ is reached, the $X$ or $Y$ coordinate is increased iteratively by one with probability $\theta$ or $1 - \theta$ respectively.

The set of all binomial sampling plans is denoted here by $\mathcal{S}$.

For any plan $S \in \mathcal{S}$ it is clear that the boundary point $\gamma \in B_S$ at which sampling stops is a sufficient statistic for $\theta$ \citep[see][p. 102]{lehmann1998}. For each such point $\gamma$, we can define $N_S(\gamma) = Y(\gamma) + X(\gamma)$, so that $N(\gamma)$ represents the total number of steps taken during sampling. An important characteristic of any plan then will be $\mathrm{E}(N_S)$, the expected number of iterations. This quantity will form the basis of comparisons of estimators across sample plans as discussed below.

If $N_S(\gamma) = n$ for some positive integer $n$ and all $\gamma \in B_\gamma$, we say $S$ is a fixed binomial sampling plan.
If, instead, $N_S(\gamma) < M$ for some positive integer $M$ and all $\gamma \in B_\gamma$, $S$ is a finite binomial sampling plan.

\section{Models}
\label{ss:models}
Here, we assume an infinite population consisting of individuals displaying some trait with probability $p$. Each individual can then be represented by independent random variables $\varphi_i, i = 1,2,3, \ldots,$ such that $\varphi_i \sim Bernoulli(p)$. Throughout, this parameter $p$ is the quantity we seek to estimate.
If group tests with pools of size $k$ are considered, we have the new random variable $\vartheta_{i}^{(k)} = \max\{\varphi_{i_1}, \varphi_{i_2}, \ldots, \varphi_{i_k}\} \sim Bernoulli(1 - q^k)$ where $q = 1 - p$. With this, we consider the three models described in the following sections.

\subsection{Model (a) - fixed binomial sampling}
\label{ss:modela}
Suppose we observe $\vartheta^{(k)}_1, \vartheta^{(k)}_r, \ldots, \vartheta^{(k)}_n$ and define $X = \sum_{i=1}^n \vartheta_i^{(k)}$. Then, $X \sim Bin(n, 1 - q^k)$ so that
\[P\left(X = x\right) = \binom{n}{x}(1 - q^k)^x(q^k)^{n - x}, \ x = 0, 1, \ldots, n.\]
In terms of the above notation, this is equivalent to the fixed binomial sampling plan $S_a$ with $\theta_{S_a} = 1 - q^k$, $B_{S_a} = \{\gamma: Y(\gamma) + X(\gamma) = n\}$ and $\mathrm{E}(N_{S_a}) = n$.

\subsection{Model (b) - inverse binomial sampling (positive)}
\label{ss:modelb}
Alternatively, if we observe $\vartheta^{(k)}_1, \vartheta^{(k)}_2, \ldots$, until the $c$th positive we can define $Y$ to be the number of negative groups until this event occurs. As such, $Y \sim NB(c, 1 - q^k)$ and so
\[P\left(Y = y\right) = \binom{c + y - 1}{y}(1 - q^k)^c(q^k)^{y},\ y = 0, 1, 2, \ldots\]
This is the same as a binomial sampling plan $S_b$ with  $\theta_{S_b} = 1 - q^k$, $B_{S_b}~=~\{\gamma: X(\gamma) = c\},$ and $\mathrm{E}(N_{S_b}) = c + \frac{cq^k}{1 - q^k} = \frac{c}{1 - q^k}$.

\subsection{Model (c) - inverse binomial sampling (negative)}
\label{ss:modelc}
If instead we observe $\vartheta^{(k)}_1, \vartheta^{(k)}_2, \ldots$, until the $c$th negative we can define $Z$ to be the number of positive groups prior to this event. As such, $Z \sim NB(c, q^k)$ and so
\[P\left(Z = z\right) = \binom{c + z - 1}{z}(q^k)^c(1 - q^k)^{z},\ z = 0, 1, 2, \ldots\]
Again, this can be expressed as a binomial sampling plan $S_c$ with  $\theta_{S_c}=1-q^k$, $B_{S_c}=\{\gamma:Y(\gamma)=c\},$ and $\mathrm{E}(N_{S_c}) = c + \frac{c(1 - q^k)}{q^k} = \frac{c}{q^k}$.

We note that, while all the previous work on group testing with inverse binomial sampling has utilized model (b) \citep[see][]{ku2006, pritchard2011a, hepworth2013}, there are many realistic cases in which model (c) would be preferable.  Moreover, we will show, based on results in \citet{degroot1959}, that an unbiased estimator for $p$ does not exist under model $(b)$ but can be constructed  under $(c)$.

Since we will present groups of estimators which apply to several of the  above sampling plans, the model used for a specific function will always be denoted in a subscript. For example, $\hat{p}_{(a)}$ will denote an estimator under model $a$.

For all models and estimators, the main quantities of interest will be the bias and MSE, given as follows,
for an estimator $\hat{p}(x)$ where $x$ is a random variable with support $\mathcal{X}$ and likelihood $L(x, p)$,
\begin{equation}
\label{eqn:bias}
\mathrm{bias}_{p}(\hat{p}) = \sum_{x \in \mathcal{X}} (\hat{p}(x) - p)L(x, p),
\end{equation}
\begin{equation}
\label{eqn:mse}
\mathrm{MSE}_{p}(\hat{p}) = \sum_{x \in \mathcal{X}} (\hat{p}(x) - p)^2L(x, p).
\end{equation}

\section{Unbiased estimation for the group testing problem}
\label{s:fsp}
As mentioned above, it is known that with standard binomial sampling no unbiased estimator exists in the group testing problem. This can be extended to all finite binomial sampling plans based on the following lemma whose proof is in Appendix \ref{app:lem1}.
\begin{lemma}
	\label{lm:1}
	Let $\mathcal{F}$ be the set of all finite binomial sampling plans with probability of success $\theta$, and $k$ any positive integer greater than one. Then, there does not exist an estimator $f$ under any sampling plan $F \in \mathcal{F}$ such that $f$ is an unbiased estimator of $\theta^{1 / k}$ or $(1 - \theta)^{1 / k}$.
\end{lemma}

For the group testing problem, where $\theta = 1 - (1-p)^k$ or $\theta = (1 - p)^k$, it follows immediately that the non-existence of an unbiased estimator of $p$ extends to this broader class of sampling plans as well.

We add here that unbiased estimation under binomial sampling plans has been studied for a range of different functions. A review can be found in \citet{sinha1992}. An alternative approach to estimating a function of the form $p^\alpha, \alpha > 0$ (of which the group testing problem is a special case) using a randomized binomial sampling scheme is presented in \citet{banerjee1979}.

\section{Estimators under model (a)}
\label{s:fbe}
As mentioned above, the MLE under model (a) has been in use for a long time and is popular, despite its often large bias, due to its simplicity and good asymptotic properties. The estimator as given in \citet{gibbs1960} is
\[\hat{p}_{MLE(a)}(x) = 1 - \left(1 - \frac{x}{n}\right)^{1/k}.\]

An alternative estimator was proposed by \citet{burrows1987} which reduces the bias by eliminating terms of $O(1 / n)$ in the Taylor expansion of the expectation of the MLE. The resulting estimator is as follows
\[\hat{p}_{B(a)}(x) = 1 - \left(1 - \frac{x}{n + \nu}\right)^{1 / k},\ \nu = \frac{k-1}{2k}.\]

As pointed out in the introduction, Burrows was able to show empirically for a range of $p$ between $0.01$ and $0.5$ and $n$ between $10$ and $200$ that, in addition to reduced bias, his estimator performed better than the MLE in terms of MSE.
The magnitude of the differences observed in his original work raise questions about the suitability of the MLE, particularly for small sample sizes. For example, with $n=10$ and $p=0.1$, the MSE for the Burrows' estimator is only one fifth that of the MLE. For larger samples, this is mitigated by the strong asymptotic performance of the MLE, although the Burrows estimator did at least slightly better in all cases considered.

\section{Estimators under models (b) and (c)}
\label{s:ibe}
The MLE under model (b), as shown in \citet{katholi2006}, is
\[\hat{p}_{MLE(b)}(y) = 1 - \left(1 - \frac{c}{y + c}\right)^{1 / k},\] and it can similarly be shown that the MLE under model (c) is

\[\hat{p}_{MLE(c)}(z) = 1 - \left(\frac{c}{z + c}\right)^{1 / k}.\]

Similar to the fixed sampling case, both MLEs presented here are biased and several alternatives have been given for model (b) to address this issue as discussed in the following sections.

\subsection{Pritchard and Tebbs estimators}
\label{ss:pt}
\citet{pritchard2011a} suggested the following three shrinkage-type estimators under model (b),
\[\hat{p}_{\alpha(b)}(y) = 1 - \left[1 - \alpha\left(\frac{c}{y + c}\right)\right]^{1 / k},\ 0 \leq \alpha \leq 1,\]
\[\hat{p}_{\beta(b)}(y) = 1 - \left[1 -\frac{c + 1}{y + c + \beta}\right]^{1 / k},\ 1 \leq \beta,\] and
\[\hat{p}_{C(b)}(y) = 1 - \left[1 -\alpha_C\left(\frac{c + 1}{y + c + \beta_C}\right)\right]^{1 / k},\ 0 \leq \alpha_C \leq 1, \ 1 \leq \beta_C,\] where each of $\alpha, \beta, \alpha_C,$ and $\beta_C$ are found numerically as the values which minimize the MSE for each respective estimator.

While the conceptual basis for these estimators is quite simple, in practice, since the MSE is a function of both the parameter and a given estimator (the expression can be found in (\ref{eqn:mse})), the optimization required to compute them depends on prior knowledge of $p$. This means the performance of the estimators can vary widely depending on the availability and accuracy of such information. In the original work, Pritchard and Tebbs recommended an upper bound $p_0$ be selected based on which the optimization should be done. This is the approach taken here where we consider estimators based on  $p_0 = .01, .1, $ and $.5$  and compare them across the full range of $p$. This allows for a basic understanding of the impact incomplete knowledge of $p$ has on the estimators. To make it clear which value of $p_0$ is being used for a given estimator, we place it in the superscript so that $\hat{p}_{C(b)}^{(.01)}$ is the estimator for method $C$ under model (b) with $p_0 = .01$.

It was pointed out by the original authors that $\hat{p}_{\beta}$ and $\hat{p}_{C}$ both outperformed $\hat{p}_{\alpha}$ in their study, with a preference for $\hat{p}_{C}$ for some choices of $k$ and $c$. As a result, to make comparisons among estimators more succinct, we only include $\hat{p}_C$ in our calculations.

\textit{Modification for Model (c)}
\label{sss:ptc}
We can modify the idea of Pritchard and Tebbs for model (c) yielding the following three estimators,
\[\hat{p}_{\alpha(c)}(z) = 1 - \left[\alpha\left(\frac{c}{z + c}\right)\right]^{1 / k},\ 0 \leq \alpha \leq 1,\]

\[\hat{p}_{\beta(c)}(z) = 1 - \left[\frac{c + 1}{z + c + \beta}\right]^{1 / k},\ 1 \leq \beta,\] and
\[\hat{p}_{C(c)}(z) = 1 - \left[\alpha_C\left(\frac{c + 1}{z + c + \beta_C}\right)\right]^{1 / k},\ 0 \leq \alpha_C \leq 1, \ 1 \leq \beta_C.\]

\subsection{Hepworth estimator}
\label{ss:hep}
Another alternative was proposed by \citet{hepworth2013} based on the Gart bias correction introduced in \citet{gart1991}. The result due to Gart states that, for a single parameter model based on independent observations, the bias of the MLE, excluding terms of $O(1 / n^2)$ (in the sequential case $O(1 / (\mathrm{E}[N])^2)$, is given by
\[B(p) = -\frac{2 \frac{dI}{dp} + \mathrm{E}\left[\frac{d^3\ell}{dp^3}\right]}{2(I(p))^2},\]
where $I(p)$ and $\ell$ are the Fisher information (contained in the sample) and log-likelihood respectively. Hepworth showed that, for model (b),
\begin{align*}
&I(p) = \frac{ck^2q^{k - 2}}{(1 - q^k)^2}, \\
&\frac{dI}{dp} = -ck^2\left(\frac{(k - 2)q^{k-3} + (k + 2)q^{2k - 3}}{(1 - q^k)^3}\right), \\
\text{and} \\
&\mathrm{E}\left[\frac{d^3\ell}{dp^3}\right] = \frac{ck}{q^3}\left(\frac{k(k + 1)q^k(1 - q^k) + 2(kq^k + q^k - 1)^2}{(1 - q^k)^3} - \frac{2}{1 - q^k}\right).
\end{align*}
He then suggests using the plugin estimator
\[\hat{p}_{G(b)}(y) = \hat{p}_{MLE(b)}(y) - B(\hat{p}_{MLE(b)}(y)).\]
Since this estimator is not defined when $\hat{p}_{MLE(b)}(y) = 1$, which occurs when $y = 0$, the value $\hat{p}_{G(b)}(0) = 1 - \left(\frac{k - 1}{2kc + k - 1}\right)^{1 / k}$ is used, which is equivalent to the Burrows correction in the fixed binomial case.

Of course, since a plugin estimator is used, this does not remove all of the $O(1 / \mathrm{E}[N])$ terms from the bias, but it does lead to a significant reduction. In the next section, we will construct an estimator for which this bias is removed analytically.

\textit{Modification for Model (c)}
\label{sss:hepc}
For model (c) we can similarly show that
\begin{align*}
&I(p) = \frac{ck^2}{q^2(1 - q^k)}, \\
&\frac{dI}{dp} = ck^2\left(\frac{2 - (2 + k)q^k}{q^3(1 - q^k)^2}\right), \\
\text{and} \\
&\mathrm{E}\left[\frac{d^3\ell}{dp^3}\right] = \frac{kc}{q^3}\left(\frac{2k^2q^{2k}}{(1 - q^k)^2} + \frac{3k(k - 1)q^k}{(1 - q^k)}  + k(k-3)\right).
\end{align*}
The estimator for this model, $\hat{p}_{G(c)}(z)$, can then be defined exactly as above with
\[\hat{p}_{G(c)}(y) = \hat{p}_{MLE(c)}(y) - B(\hat{p}_{MLE(c)}(y)),\] using the above values in the definition of $B(p)$.
In this case, the bias correction is not defined when $p = 0$, so that the value $\hat{p}_{G(c)}(0) = 0$, which is also equivalent to the Burrows estimator in the fixed case, is used here.

\subsection{Extension of Burrows to models (b) and (c)}
\label{ss:bur}
In this section, we extend the idea of Burrows in the fixed sampling case to the sequential models discussed here, with the modification that we seek to remove terms of order $O(1 / \mathrm{E}[N])$ from the bias. This is in the same spirit as the Gart estimator presented by Hepworth, but it has the added advantage of producing estimators which achieve the desired reduction in bias theoretically, not only as an approximation.

\subsubsection{Burrows type estimator for model (b)} 
To apply the idea of Burrows to model (b) we begin with the modified MLE,
\[\hat{p}_{B(b)} = 1 - \left(\frac{y + \nu}{y + c + \eta}\right)^{1 / k}.\] Then, we use the Taylor expansion for the expectation of $\hat{p}_{B(b)}$ to find all terms of $O(1 / \mathrm{E}[N])$, and solve for $\nu$ and $\eta$ which result in the removal of such terms. The result is $\nu = \eta + 1 = \frac{k - 1}{2k}$, yielding the estimator
\begin{equation}
\label{eqn:burb}
\hat{p}_{B(b)} = 1 - \left(\frac{y + \nu}{y + c + \nu - 1}\right)^{1 / k},\ \nu = \frac{k - 1}{2k}.
\end{equation}
A formal proof is provided in Appendix \ref{app:burb}.

It should be noted that if $c = 1$, this estimator will trivially yield zero for all values of $y$ making it unusable in such a case. For group testing problems, when $p$ is generally small, this is unlikely to be an issue in applications as the choice of $c=1$ will often yield unreasonably small expected sample sizes. For example, if $p=0.1$ and $k = 2,$ then $c=1$ yields $\mathrm{E}[N] = \frac{100}{19},$ with this value decreasing for all $k > 2$.

\subsubsection{Burrows type estimator for model (c)} 
The application to model (c) follows exactly as above beginning with the estimator
\[\hat{p}_{B(c)} = 1 - \left(\frac{c + \nu}{y + c + \eta}\right)^{1 / k}.\] The desired solution is $\nu = \eta = \frac{k-1}{2k} - 1$ so the estimator is
\begin{equation}
\label{eqn:burc}
\hat{p}_{B(c)} = 1 - \left(\frac{c +\nu - 1}{y + c + \nu -1}\right)^{1 / k},\ \nu = \frac{k-1}{2k}.
\end{equation}
The proof of (\ref{eqn:burc}) is nearly identical to that of (\ref{eqn:burb}) and is omitted.
\subsection{Construction of unbiased estimator}
\label{ss:deg}
In this section, we show how to construct an unbiased estimator for $p$ under model (c) directly from a theorem due to Degroot. The ability to apply this result to a similar function was mentioned by \citet{hall1963}. The result is given in the following Theorem and can be found, along with its proof, in Theorem 4.1 of \citet{degroot1959}.
\begin{theorem}[Degroot]
	\label{thm:deg}
	Let $W \sim NB(c, 1 - \theta)$. Then, a function $h(\theta)$ is estimable unbiasedly if and only if it can be expanded in a Taylor series on the interval $|\theta| < 1$. If $h(\theta)$ is estimable unbiasedly, then its unique estimator is given by
	\[\hat{h}(w) = \frac{(c - 1)!}{(w + c - 1)!}\frac{d^w}{d\theta^w}\left[\frac{h(\theta)}{(1 - \theta)^c}\right]_{\theta = 0}, \ w = 0, 1, 2, \ldots.\]
\end{theorem}

To apply this under model (c) we have $\theta_c = 1 - q^k$ and want to estimate $h_c(\theta) = (1 - \theta)^{1 / k} = q$ which leads immediately to an estimate of $p$. Then Theorem \ref{thm:deg} yields the estimator
\[\hat{p}_{D(c)}(z) = \left\{\begin{array}{ll} 0, & z = 0, \\ 1 - \prod_{j=1}^z \left(\frac{j + c - 1 - 1 / k}{j + c - 1}\right), & z = 1, 2, 3, \ldots.\end{array}\right.\]

It should be noted that, under model (b), we have $\theta_b = q^k$ so that $h_b(\theta) = \theta^{1 / k} = q$. However, $h_b$ does not have a Taylor expansion at the point $\theta = 0$, so by Degroot's Theorem no unbiased estimator exists under this model.

\section{Numerical comparisons}
Comparisons among the estimators presented here can be challenging due to the number of variables which must be considered. These include $p$, $\mathrm{E}(N)$ (which reduces to the fixed $n$ under model (a)), and $k$.  This is true even among estimators in the same model where, for example, with $p$ and $\mathrm{E}(N)$ fixed, one estimator may perform much better for one choice of $k$ while a second estimator is superior with a different selection. As such, any predetermined values may unduly favor one estimator over another.

To deal with this, for all comparisons we considered $p$ and $\mathrm{E}(N)$ fixed and then chose the value of $k$ for each estimator which yields the smallest MSE. This was done using a basic grid search over $k =\{2,\ldots, 50\}$. The choice of 50 as an upper bound is somewhat arbitrary, but reflects the vast majority of group testing applications in which even smaller pool sizes are used. It should be noted that, particularly for the Burrows and Gart estimators, it may be possible to choose a value of $k$ which has a much smaller bias without inflating the MSE when compared to the minimizing value. This is especially an issue for small $\mathrm{E}(N)$ and is one more design issue which must be considered in application.

For the other values, we looked at $p = 0.01, 0.05, 0.1, 0.2, 0.3,$ and $0.5$ as well as $\mathrm{E}[N]= 25$ and $100$.

Note that the value of $c$ for models (b) and (c) is completely determined by $\mathrm{E}(N)$, $p$, and $k$ as shown in sections \ref{ss:modelb} and \ref{ss:modelc}. However, since $c$ is constrained to the positive integers, we always selected $c$ to be the largest integer such that $\mathrm{E}_c(N) \leq \mathrm{E}(N)$ where the latter is the original target value. Since the sequential plans will rarely have a $c$ value which allows them to attain the upper bound exactly, this approach will yield a slight advantage to estimators under the fixed plans.

For models (b) and (c), since neither of equations (\ref{eqn:bias}) or (\ref{eqn:mse}) can  be calculated exactly, we instead found $\nu_b$ such that $P(Y > \nu_b) \leq 1\times10^{-6}$ and $\nu_c$ such that $P(Z > \nu_b) \leq 1\times10^{-6}$ and took these values as the upper limits in the appropriate sums.

Tables \ref{tab:RB25} and \ref{tab:RB100} contain comparisons for the relative bias defined, for an estimator $\hat{p}$, as $100 \times \frac{\mathrm{E}[\hat{p} - p]}{p}$.

\begin{table}[ht]
	\caption{Relative bias comparisons for $\mathrm{E}(N) = 25$.}
	\begin{tabular*}{\textwidth}{l@{\extracolsep{\fill}}*{7}{S[table-format=-3.4]}}
		\toprule
		$\hat{p}\backslash p$& 0.01 & 0.05 & 0.1 & 0.2 & 0.3 & 0.5 \\
		\midrule
		$\hat{p}_{MLE(a)}$       & 2.6656 & 3.2943 & 3.0528 & 2.8887 & 2.1418 & 1.6494 \\
		$\hat{p}_{MLE(b)}$       &14.2520 & 9.9336 & 8.6084 & 7.3317 & 5.8485 & 4.0555 \\
		$\hat{p}_{MLE(c)}$       &-2.5988 & -5.0153 & -4.1250 & -3.9986 & -4.6445 & -4.7977 \\
		$\hat{p}_{B(a)}$         & 0.0240 & 0.0930 & -10.8749 & -11.1933 & -9.1894 & -8.9490 \\
		$\hat{p}_{B(b)}$         & 0.0466 & -1.4009 & -10.8632 & -9.9224 & -10.2183 & -9.8681 \\
		$\hat{p}_{B(c)}$         & 0.0241 & 0.1182 & 0.0666 & 0.0506 & 0.0721 & 0.0513 \\
		$\hat{p}_{C(b)}^{(.01)}$ & -12.2992 & -13.1649 & -13.3573 & -11.2163 & -15.5832 & -31.6441 \\
		$\hat{p}_{C(b)}^{(.1)}$  & 21.0588 & 2.0035 & -4.2050 & -13.7477 & -31.5621 & -44.0850 \\
		$\hat{p}_{C(b)}^{(.5)}$  & 22.4240 & 12.0359 & 22.3028 & -3.6826 & -3.1759 & -3.2645 \\
		$\hat{p}_{C(c)}^{(.01)}$ & 0.0695 & -8.2602 & -9.9764 & -9.5557 & -10.9204 & -10.3177 \\
		$\hat{p}_{C(c)}^{(.1)}$  & 19.9813 & 4.5841 & 3.1493 & -6.9119 & -8.0160 & -8.9280 \\
		$\hat{p}_{C(c)}^{(.5)}$  & 172.1298 & 34.3271 & 118.5546 & 48.1304 & 24.9428 & 2.7856 \\
		$\hat{p}_{G(b)}$         & -0.0975 & -2.4344 & -12.1228 & -12.2914 & -10.1906 & -9.6814 \\
		$\hat{p}_{G(c)}$         & -0.0467 & -0.2033 & -0.1358 & -0.1345 & -0.1963 & -0.2324 \\
		\bottomrule
	\end{tabular*}

	\label{tab:RB25}
\end{table}

\begin{table}[ht]
	\caption{Relative bias comparisons for $\mathrm{E}(N) = 100$.}
	\begin{tabular*}{\textwidth}{l@{\extracolsep{\fill}}*{7}{S[table-format=-3.4]}}
		\toprule
		$\hat{p}\backslash p$& 0.01 & 0.05 & 0.1 & 0.2 & 0.3 & 0.5 \\
		\midrule
		$\hat{p}_{MLE(a)}$       & 0.6411 & 1.0847 & 1.0412 & 0.9556 & 0.7111 & 0.3821 \\
		$\hat{p}_{MLE(b)}$       & 2.7318 & 1.8045 & 1.7197 & 1.5546 & 1.3314 & 0.8894 \\
		$\hat{p}_{MLE(c)}$       & -0.6610 & -1.0926 & -1.1110 & -1.1038 & -1.1573 & -1.1328 \\
		$\hat{p}_{B(a)}$         & 0.0014 & 0.0059 & 0.0061 & 0.0054 & 0.0039 & 0.0025 \\
		$\hat{p}_{B(b)}$         & 0.0017 & 0.0055 & 0.0056 & 0.0059 & 0.0076 & 0.0026 \\
		$\hat{p}_{B(c)}$         & 0.0015 & 0.0046 & 0.0045 & 0.0039 & 0.0039 & 0.0025 \\
		$\hat{p}_{C(b)}^{(.01)}$ & -8.6125 & -10.8966 & -10.8605 & -9.7154 & -9.6308 & -13.0284 \\
		$\hat{p}_{C(b)}^{(.1)}$  & 4.0879 & -1.3372 & -2.9799 & -3.5816 & -4.3240 & -3.7088 \\
		$\hat{p}_{C(b)}^{(.5)}$  & 4.6278 & 2.2599 & 2.2304 & -5.1475 & -4.0225 & -1.6152 \\
		$\hat{p}_{C(c)}^{(.01)}$ & 0.0304 & -2.1660 & -2.6021 & -2.6957 & -2.7848 & -2.5173 \\
		$\hat{p}_{C(c)}^{(.1)}$  & 3.5657 & 0.9416 & 2.2391 & -1.3774 & -2.1000 & -2.2255 \\
		$\hat{p}_{C(c)}^{(.5)}$  & 28.5418 & 5.6584 & 2.7947 & 1.4592 & 1.3026 & 1.8734 \\
		$\hat{p}_{G(b)}$         & -0.0027 & -0.0102 & -0.0075 & -0.0042 & -0.0029 & 0.0053 \\
		$\hat{p}_{G(c)}$         & -0.0030 & -0.0094 & -0.0100 & -0.0104 & -0.0121 & -0.0127 \\
		\bottomrule
	\end{tabular*}
	\label{tab:RB100}
\end{table}

From Table \ref{tab:RB25} we see that, for $\mathrm{E}(N) = 25$, the Burrows and Gart estimators significantly reduce the bias under each model when $p$ is small, and always under model (c). The increase in bias under models (a) and (b) for larger $p$ can be avoided, as mentioned above, by carefully selecting an alternative value of $k$ which reduces the bias with only a small increase in the MSE. We do not attempt this in the tables here so as to maintain consistency across comparisons. It should be noted that this is generally not possible for the other estimators considered. As expected, the Burrows estimators do have smaller bias (generally a reduction of half or more) when compared with the Gart estimators for small $n$. For the estimator due to Pritchard and Tebbs, we see that, while the bias is well controlled when the true $p$ is near the known upper bound $p_0$, as $p$ moves away from this value it tends to become, often significantly, inflated.

In Table \ref{tab:RB100}, with $\mathrm{E}(N) = 100$,  we see the same basic trends as in the previous Table, but with the expected decrease in bias across all estimators. One significant difference we see is that, for some larger values of $p$, the Gart estimator now has smaller absolute bias under model (b) than the Burrows estimator. This will likely continue to occur as $n$ increases, although the magnitude of the bias in such cases will be significantly small.

It is also interesting to note that, especially for the smaller $\mathrm{E}(N)$, the bias tends to be much smaller for estimators under model (c) relative to their counterparts under model (b).

Comparisons based on MSE can be found in Tables \ref{tab:MSE25} and \ref{tab:MSE100}, where each value is multiplied by $10,000$ for ease of interpretation.

\begin{table}[ht]
	\caption{MSE comparisons for $\mathrm{E}(N) = 25$ ($10000 \times MSE)$.}
	\begin{tabular*}{\textwidth}{l@{\extracolsep{\fill}}*{7}{S[table-format=-3.4]}}
		\toprule
		$\hat{p}\backslash p$& 0.01 & 0.05 & 0.1 & 0.2 & 0.3 & 0.5 \\
		\midrule
		$\hat{p}_{MLE(a)}$       & 0.1119 & 1.9982 & 7.3243 & 24.5901 & 46.1621 & 82.4696 \\
		$\hat{p}_{MLE(b)}$       & 1.3059 & 4.9489 & 12.8547 & 38.6643 & 62.1209 & 101.5301 \\
		$\hat{p}_{MLE(c)}$       & 0.1010 & 1.6105 & 6.0341 & 22.6446 & 43.7033 & 96.6345 \\
		$\hat{p}_{B(a)}$         & 0.1039 & 1.6010 & 3.6165 & 13.2301 & 26.7432 & 56.3798 \\
		$\hat{p}_{B(b)}$         & 0.1477 & 1.5911 & 4.8515 & 17.2066 & 33.6451 & 64.2978 \\
		$\hat{p}_{B(c)}$         & 0.1046 & 1.6237 & 6.1142 & 22.8252 & 42.7642 & 90.0256 \\
		$\hat{p}_{C(b)}^{(.01)}$ & 0.0668 & 1.2629 & 4.8836 & 17.6856 & 39.9143 & 258.7384 \\
		$\hat{p}_{C(b)}^{(.1)}$  & 0.2604 & 1.5266 & 1.8650 & 18.2526 & 94.5432 & 489.9991 \\
		$\hat{p}_{C(b)}^{(.5)}$  & 0.2840 & 3.3462 & 9.8075 & 2.3389 & 3.1281 & 12.9734 \\
		$\hat{p}_{C(c)}^{(.01)}$ & 0.0925 & 1.5588 & 6.2914 & 24.0942 & 50.3361 & 116.7416 \\
		$\hat{p}_{C(c)}^{(.1)}$  & 0.1796 & 1.8020 & 1.1595 & 18.0997 & 40.8000 & 103.1921 \\
		$\hat{p}_{C(c)}^{(.5)}$  & 3.8199 & 7.4057 & 147.6589 & 106.5589 & 76.7909 & 17.5551 \\
		$\hat{p}_{G(b)}$         & 0.1455 & 1.6526 & 5.3620 & 18.9407 & 36.8691 & 69.4281 \\
		$\hat{p}_{G(c)}$         & 0.1045 & 1.6219 & 6.1097 & 22.8195 & 42.8159 & 90.4985 \\
		$\hat{p}_{D(c)}$         & 0.1046 & 1.6230 & 6.1124 & 22.8217 & 42.7695 & 90.0741 \\
		\bottomrule
	\end{tabular*}
	\label{tab:MSE25}
\end{table}

\begin{table}[ht]
	\caption{MSE comparisons for $\mathrm{E}(N) = 100$ ($10000 \times MSE)$.}
	\begin{tabular*}{\textwidth}{l@{\extracolsep{\fill}}*{7}{S[table-format=-3.4]}}
		\toprule
		$\hat{p}\backslash p$& 0.01 & 0.05 & 0.1 & 0.2 & 0.3 & 0.5 \\
		\midrule
		$\hat{p}_{MLE(a)}$       & 0.0261 & 0.3939 & 1.4850 & 5.2302 & 10.0945 & 19.1070 \\
		$\hat{p}_{MLE(b)}$       & 0.0298 & 0.4183 & 1.5700 & 5.4881 & 10.6862 & 19.6360 \\
		$\hat{p}_{MLE(c)}$       & 0.0257 & 0.3766 & 1.4540 & 5.1569 & 10.0262 & 19.8416 \\
		$\hat{p}_{B(a)}$         & 0.0257 & 0.3737 & 1.4140 & 4.9961 & 9.7855 & 18.8259 \\
		$\hat{p}_{B(b)}$         & 0.0274 & 0.3862 & 1.4547 & 5.1179 & 10.1105 & 19.1499 \\
		$\hat{p}_{B(c)}$         & 0.0259 & 0.3772 & 1.4547 & 5.1431 & 9.9391 & 19.3963 \\
		$\hat{p}_{C(b)}^{(.01)}$ & 0.0285 & 0.6530 & 2.5048 & 9.0790 & 18.1887 & 54.2248 \\
		$\hat{p}_{C(b)}^{(.1)}$  & 0.0311 & 0.3584 & 1.2747 & 4.7773 & 9.7010 & 19.6961 \\
		$\hat{p}_{C(b)}^{(.5)}$  & 0.0320 & 0.4214 & 1.5884 & 3.8762 & 5.4868 & 17.4286 \\
		$\hat{p}_{C(c)}^{(.01)}$ & 0.0251 & 0.3780 & 1.4839 & 5.3324 & 10.5320 & 21.2195 \\
		$\hat{p}_{C(c)}^{(.1)}$  & 0.0295 & 0.3673 & 0.5823 & 4.8543 & 9.8917 & 20.4714 \\
		$\hat{p}_{C(c)}^{(.5)}$  & 0.2448 & 0.9582 & 2.0042 & 5.0234 & 9.6493 & 8.4774 \\
		$\hat{p}_{G(b)}$         & 0.0274 & 0.3857 & 1.4528 & 5.1122 & 10.0963 & 19.1494 \\
		$\hat{p}_{G(c)}$         & 0.0259 & 0.3772 & 1.4546 & 5.1433 & 9.9402 & 19.4030 \\
		$\hat{p}_{D(c)}$         & 0.0259 & 0.3772 & 1.4546 & 5.1431 & 9.9392 & 19.3969 \\
		\bottomrule
	\end{tabular*}

	\label{tab:MSE100}
\end{table}

Both tables, based on $\mathrm{E}(N)= 25$ and $\mathrm{E}(N) = 100$ respectively, show nearly identical patterns among the estimators, with only the magnitude of the values decreasing with the expected sample size. For the Pritchard and Tebbs estimator, we see that, when $p_0$ is chosen close to the true $p$, it always outperforms the other estimators in terms of MSE. However, as $p$ moves away from $p_0$, we again see that the value becomes inflated, almost always making it the worst estimator among those considered.

For all other estimators, we see that the Burrows, Gart, and Degroot estimators generally outperform the MLE under each appropriate model. The exception to this is model (c), for which the MLE actually yields a smaller MSE when compared to the alternatives for all values of $p$ but $p=0.5$. Similar to the bias comparisons, we see that estimators under model (c) all have smaller MSE than those under model (b) when $p$ is very small, but this is reversed as $p$ increases (the exception is the MLE for which the estimator under model (c) is always better than the one under model (b)). Interestingly, the most consistently small MSE is achieved by the Burrows estimator under model (a). This, however, should be considered in light of the above mentioned advantage for estimators under the fixed sample plans due to the slightly larger expected sample sizes.

\section{Discussion} 
\label{s:dis}
The primary results of this work can be divided into two parts. First, we have shown that, despite decades of attempts at bias minimization, an unbiased estimator does exist for the group testing problem if one is willing to consider inverse binomial sampling based on counting negatives. Our numerical comparisons indicate that this is achieved together with a significant reduction in the MSE relative to the MLE for the standard fixed binomial case, as well as a comparable reduction to the alternative estimators found in the literature.

Likewise, we have provided proofs that, under the two most prominent models in the literature, models (a) and (b), there exist no unbiased estimators. As such, it follows that the Degroot estimator introduced here is trivially the optimal unbiased estimator across the class of sampling plans considered. These facts combined make this estimator particularly desirable, at least from a theoretical standpoint.

That said, if some bias can be tolerated, and there is complete freedom to choose a model, there seems to be little reason to suggest moving away from fixed binomial sampling. Not only is this model much more familiar, but the Burrows estimator outperforms in terms of MSE, at least slightly, all other estimators across sampling plans in our comparisons.

The second area relates to situations in which inverse binomial sampling is desired a priori. \citet{pritchard2011a} give some motivation for problems in which testing groups until $c$ positives are attained, which corresponds to our model $(b)$. For this situation, we have extended the idea of Burrows to provide an estimator which removes terms of $O(1 / \mathrm{E}(N))$ from the bias while also yielding drastic reductions in the MSE when compared with the MLE. While this is similar in motivation to the estimator proposed in \citet{hepworth2013}, the form of the estimator is much simpler and has the advantage of achieving the desired bias reduction analytically. In our comparisons, both estimators performed very similarly, with the Burrows type estimator generally doing better in terms of bias. In contrast, our results indicate that the estimators suggested in \citet{pritchard2011a}, while yielding extremely small MSE values in some cases, are, in addition to being heavily biased, too dependent on precise prior knowledge of $p$ to be useful in most applications involving small sample sizes.

When comparing models (b) and (c), it is clear that neither dominates the other uniformly in the considered comparisons. In fact, we see that in each case estimators under one model slightly outperform their counterparts in exactly one of bias or MSE (the exception being the MLE for which model (c) is always better when $p$ is not too large). Furthermore, the differences tend to be very small and in application, when $p$ is unknown, it would likely be impossible to determine which is optimal for any given set of design constraints.
As such, if an application exists for which sampling until $c$ negatives is possible, we believe that model (c) with the Degroot estimator will usually be preferred. Aside from the desirable theoretical properties mentioned above, the unbiasedness of this estimator can greatly simplify the design process by allowing one to focus solely on the MSE.

We emphasize that the perspective of this work has been one of design, in the sense of choosing a study design which yields optimal results in terms of reduced MSE and bias (as a secondary goal), particularly when sample sizes are relatively small. This is in contrast with the problem, often encountered in applications, of selecting a best estimator when the study design is constrained by experimental or other factors. As such, while the results presented here do give a general indication of the small sample performance of each estimator, any actual application will likely require a more targeted comparison to choose the best design and analysis strategy for the specific question at hand.

Potential future work in sequential group testing remains very broad, since the area is largely undeveloped. In particular, questions of design such as the optimal choice of group size, possibly using an adaptive approach, are extremely important if such designs are to be utilized in applications. Furthermore, extensions to further cases such as the presence of misclassification and multiple group sizes in small sample studies are vital. Additionally, there is no reason to limit research to the inverse binomial model, and generalizations to other sequential sampling plans could be explored, although such work would likely be very application specific. Of more theoretical interest, the possibility of determining all possible designs yielding an unbiased estimator for $\theta^{1/k}$, where $\theta$ is a Bernoulli parameter and $k$ is a positive integer greater than one, as was done for the case $1 / \theta$ \citep[see][]{sinha1985} remains open and appears to be a difficult problem.

\section*{Acknowledgement}
We thank Yan Lumelskii, Aiyi Liu, and Graham Hepworth for useful discussions as well as the associate editor and editor in chief for their time and advice.
We also thank Shelemyahu Zacks and Albert Vexler for helpful comments on the manuscript.
The work of the third author was supported by the National Cancer Institute Intramural Program.

\appendix
\section{Proofs}
\subsection{Proof of Lemma \ref{lm:1}}
\label{app:lem1}
Let $F$ be a finite sequential binomial sampling plan with probability of success $\theta$ and define $K(\gamma)$ to be the number of ways to reach the boundary point $\gamma \in B_S$. Since $B_S$ is finite, set $\eta = \max_{\gamma \in B_S}\{Y(\gamma) + X(\gamma)\}$ and we have, for any statistic $T$ defined on $B_S$,
\[\mathrm{E}(T) = \sum_{\gamma \in B_S} T(\gamma)K(\gamma) \theta^{X(\gamma)}(1 - \theta)^{Y(\gamma)} = \sum_{i=0}^\eta C_j \theta^j,\] for some constants $C_0, \ldots, C_\eta$ which is a polynomial of degree at most $\eta$. Since it is impossible for this sum to equal $\theta^{1 / k}$ or $(1 - \theta)^{1 / k}$ for all $\theta$, it follows that no unbiased estimator exists for either function.

\subsection{Proof of (\ref{eqn:burb})}
\label{app:burb}
Let $Y \sim NB(c, 1 - q^k)$. For convenience, we begin with the estimator of $q$, $\hat{q}_{B(b)}(y) = \left(\frac{y + \nu}{y + c + \eta}\right)^\xi$, where $\xi = \frac{1}{k}$ and note that $O\left(\frac{1}{\mathrm{E}(N)}\right)$ is equivalent to $O\left(\frac{1}{c}\right)$.
Then, taking the Taylor expansion of $\hat{q}_{B(b)}$ about $y_0 = \mathrm{E}(Y) = \frac{cq^k}{1 - q^k}$ yields
\begin{multline}
\label{eqn:p1}
\hat{q}_{B(b)}(y) = \left(\frac{y_0 + \nu}{y_0 + c + \eta}\right)^\xi + \xi\left(\frac{y_0 + \nu}{y_0 + c + \eta}\right)^{\xi - 1}\frac{(c + \eta - \nu)}{(y_0 + c + \eta)^2}(y - \mathrm{E}[Y])\\  +\left[\xi(\xi - 1)\left(\frac{y_0 + \nu}{y_0 + c + \eta}\right)^{\xi - 2}\frac{(c + \eta - \nu)^2}{(y_0 + c + \eta)^4}\right.\\ - \left.2\xi\left(\frac{y_0 + \nu}{y_0 + c + \eta}\right)^{\xi - 1}\frac{(c + \eta - \nu)}{(y_0 + c + \eta)^3}\right]\frac{(y - \mathrm{E}[Y])^2}{2},
\end{multline} plus terms which, after taking the expectation, will be $O\left(\frac{1}{c^2}\right)$.
Then, taking the expectation of (\ref{eqn:p1}) yields
\begin{multline}
\label{eqn:p2}
\mathrm{E}[\hat{q}_{B(b)}(Y)] = \left(\frac{y_0 + \nu}{y_0 + c + \eta}\right)^\xi+ \left[\xi(\xi - 1)\left(\frac{y_0 + \nu}{y_0 + c + \eta}\right)^{\xi - 2}\frac{(c + \eta - \nu)^2}{(y_0 + c + \eta)^4}\right. \\ \left. - 2\xi\left(\frac{y_0 + \nu}{y_0 + c + \eta}\right)^{\xi - 1}\frac{(c + \eta - \nu)}{(y_0 + c + \eta)^3}\right]\frac{cq^k}{2(1 - q^k)^2} + O\left(\frac{1}{c^2}\right).
\end{multline}
Now, taking again the Taylor expansion of (\ref{eqn:p2}) about $(\nu, \eta) = (0,0)$ and plugging in the value for $y_0$ yields
\begin{multline}
\label{eqn:p3}
\mathrm{E}[\hat{q}_{B(b)}(Y)] = q + \frac{\xi q (1 - q^k) \nu}{cq^k}  - \frac{\xi q (1 - q^k)q^k \eta}{cq^k} + \frac{\xi (\xi - 1) q (1 - q^k)^2}{2cq^k}\\ -\frac{\xi q (1 - q^k)q^k}{cq^k} + O\left(\frac{1}{c^2}\right).
\end{multline}
Then, the terms $O\left(\frac{1}{c}\right)$ will disappear from (\ref{eqn:p3}) if
\[\nu - q^k\eta + \frac{(\xi - 1)(1-q^k)}{2} - q^k = 0,\] which has the unique solution
\[\nu = \eta + 1 = \frac{k-1}{2k}.\]

\end{document}